\documentclass[suppldata]{interact}
\usepackage[autostyle]{csquotes}
\usepackage{epstopdf}
\usepackage[caption=false]{subfig}

\usepackage[numbers,sort&compress]{natbib}
\bibpunct[, ]{[}{]}{,}{n}{,}{,}
\makeatletter
\def\NAT@def@citea{\def\@citea{\NAT@separator}}
\makeatother
\usepackage{anyfontsize}
\theoremstyle{plain}
\newtheorem{theorem}{Theorem}[section]
\newtheorem{lemma}[theorem]{Lemma}
\newtheorem{corollary}[theorem]{Corollary}
\newtheorem{proposition}[theorem]{Proposition}

\theoremstyle{definition}
\newtheorem{definition}[theorem]{Definition}
\newtheorem{remark}[theorem]{Remark}
\newtheorem{example}[theorem]{Example}

\usepackage{xcolor}

\theoremstyle{remark}

\begin{document}
	

	\title{Frame Scaling  by Graphs}
  \author{\name {Ayyanar K \thanks{Ayyanar K. Email address: ayyanark.217ma003@nitk.edu.in} P. Sam Johnson \thanks{ P. Sam Johnson(corresponding author). Email address: sam@nitk.edu.in}and A. Senthil Thilak \thanks{ A. Senthil Thilak. Email address: thilak@nitk.edu.in}} \affil{Department of Mathematical and Computational Sciences, \\
			National Institute of Technology Karnataka, Surathkal, Mangaluru 575025, India.}}
  
	\maketitle

	\begin{abstract}
     In this paper, we investigate the scalability of a given frame in $\mathbb{R}^n$ by using graphs. For each frame $\phi$ in $\mathbb{R}^n$, we associate a simple undirected graph $G(\phi)$ and use it to verify the scalability of $\phi$. We provide some necessary conditions to test the scalability of a given frame. Finally, we study the scalability of some special classes of frames by using graphs.
	\end{abstract}
	
	\begin{keywords}
		Frame graph, tight frame, scalable frame, scalable graph.
	\end{keywords}
      \begin{amscode}42C15, 05C50.\end{amscode}
    
    \section{Introduction}\label{sec1}

The development of frame theory can be traced back to the mid-20th century when mathematicians began to explore alternative methods for signal
representation beyond the traditional use of bases. In 1952, Duffin and Schaeffer introduced the concept of frame theory for Hilbert spaces through their
pioneering work in \cite{duf}. However, it was not until 1986 when frame theory
gained popularity through the revolutionary work of Daubechies, Grossman
and Meyer \cite{Dau}. Today, frame theory continues to be an expanding area of research with applications in various disciplines.

	 A finite sequence of vectors $ \phi $ in $\mathbb{R}^n$ is a \textit{frame} if and only if $ \phi $ spans $ \mathbb{R}^n $ (see Section \ref{sec2} for the definition).  In recent years, frames in finite dimensional inner product spaces have received much attention from both pure and applied mathematics, as they possess a rich redundant structure when compared to bases for vector spaces. Given  a Parseval frame $\phi=\{f_i\}_{i=1}^m$ (see Section \ref{sec2} for the definition) in $\mathbb{R}^n$,  every vector (signal) $f$ in $\mathbb{R}^n$ is transmitted as a sequence of coefficients $\{\langle f, f_i\rangle\}\}_{i=1}^m$ and it can be recovered by the following reconstruction formula :\begin{equation}\label{eqn_one}
    f=\sum\limits_{i=1}^m\langle f, f_i\rangle f_i,\,\,\,\, \text{for all}\,\, f\in \mathbb{R}^n.
\end{equation} 
 If $\phi$ is a frame (not necessarily Parseval), the reconstruction formula depends of the operators  $S_\phi^{-1}$ or $S_\phi^{-1/2}$ : 
 \begin{equation}\label{eqn_two}
    f=\sum\limits_{i=1}^m\langle f, S_\phi^{-1}f_i\rangle f_i=\sum\limits_{i=1}^m\langle f, S_\phi^{-1/2}f_i\rangle S_\phi^{-1/2}f_i,\,\,\,\, \text{for all}\,\, f\in \mathbb{R}^n,
\end{equation}
  where $S_\phi$ is the frame operator of $\phi$. In other words, if  $\phi=\{f_i\}_{i=1}^m$ is a frame in $\mathbb{R}^n$, then $\{S_\phi^{-1/2}f_i\}_{i=1}^m$ is a Parseval frame.  To get a Parseval frame from a frame using the inverse of frame operator, we encounter two important issues. Firstly, computing inverse and  square root of the frame operator are not easy in higher dimensional  spaces. Secondly, the frame $\{S_\phi^{-1/2}f_i\}_{i=1}^m$  corresponding to the inverse does not preserve the original structure of the given frame $\{f_i\}_{i=1}^m$. 
  These  issues can be resolved by modifying the frames vectors in the frame $\{f_i\}_{i=1}^m$  by (non-negative) scalings $\{a_i\}_{i=1}^m$ ($a_i\geq 0)$ so as to convert $\{f_i\}_{i=1}^m$ into a Parseval frame $\{a_if_i\}_{i=1}^m$.  The process of scaling helps to retain the basic properties of the frame vectors and has the expression which is similar to (\ref{eqn_one}):
   \begin{equation}\label{eqn_three}
  f=\sum\limits_{i=1}^m\langle f, a_if_i\rangle a_if_i,\,\,\,\, \text{for all}\,\, f\in \mathbb{R}^n.
  \end{equation}

A frame is called a scalable frame if there are non-negative constants $\{a_i\}_{i=1}^m$ making $\{f_i\}_{i=1}^m$ into a Parseval frame $\{a_if_i\}_{i=1}^m$. It is mentioned in \cite{cas23} that frame scaling is a very challenging problem in frame theory. For scalable frames, we refer the reader to \cite{ cas17, cas20, cha15, che15, cop13, dom15, gkuty13, kuty14}. In this paper, we discuss scalability of frames by analyzing graphs associated with them. 

We briefly describe the contents of the paper. Section \ref{sec2} deals with the basic definitions and results on frames and graphs, which are useful in the sequel. In Section \ref{sec3}, we provide some relations between graphs and frames. In Section \ref{sec4}, we provide some necessary conditions to examine scalability (strictly scalability) by using graphs. In the end of the section, we verify the scalability of some special classes of frames.

\section{Preliminaries}\label{sec2}
In this section, we recall some basic facts about finite frames and graphs. For more details, we refer the reader to \cite{han07, west} and references therein.

A \textit{graph} $G$ is an ordered pair $ (V(G), E(G)) $, where $V(G)$ is a non-empty set of elements called \textit{vertices} and $E(G)$ is a set of elements called \textit{edges}, where each edge is an ordered or unordered pair of vertices. In this paper, we consider only undirected simple graphs (graphs with no loops and no multiple edges). The number of vertices and edges of $ G $ are respectively referred to as the \textit{order} and \textit{size} of $ G $. Let $ u, v \in V(G) $. If $ e \in E(G) $ and $ e = (u, v) $, then $ u $ and $ v $ are said to be \textit{adjacent vertices}, denoted by $ u \sim v $. An edge with $v$ as an endpoint is said to be \textit{incident on $ v $}, and the number of edges incident on $ v $ is called the \textit{degree} of $v$, denoted by $deg(v)$. A graph where each vertex has the same degree is called a \textit{regular graph}.

The collection of all vertices adjacent to $v$ is called the open neighborhood of $v$ and is denoted by $N(v)$. If $N(v)$ contains $v$, then it is called the closed neighborhood of $v$, and it is denoted by $N[v]$.  
A set of mutually non-adjacent vertices in a graph $G$ is called an \textit{independent} set of $G$. The cardinality of a maximum independent set in $G$ is called the \textit{independent number} of $G$, denoted by $\alpha(G)$.

A \textit{walk} of a graph is an alternating sequence of vertices and edges. The number of edges in a walk is called the length of the walk. A \textit{trail} is a walk with distinct edges, and a \textit{path} is a trail with distinct vertices.  A path on $n$-vertices is denoted by $P_n$. The length of the shortest path from $u$ to $v$  is called the distance of $u$ from $v$
 and is denoted by $d(u,v)$, where $u, v\in V(G)$ and $u \neq v$. The diameter of a graph $G$ is the maximum distance between the pair of vertices, and it is denoted by $diam(G)$. A path with the same initial and terminal endpoints is called a \textit{cycle}, and a cycle on $n$-vertices is denoted by $C_n$. A graph having no cycle is \textit{acyclic}. A graph $G$ is \textit{connected} if there exists a path between every pair of vertices in $G$; otherwise, it is \textit{disconnected}. A connected graph with no cycle is called a \textit{tree}.

A graph is \textit{complete} if there is an edge between every pair of vertices, and a complete graph on $n$-vertices is denoted by $K_n$. A graph with $n$ vertices and no edges is called \textit{totally disconnected} graph and it is denoted by $\overline{K_n}$. A graph $ G $ is \textit{bipartite} if $ V(G) $ is the disjoint union of two independent sets, say $ X $ and $ Y $. Here, $ X $ and $ Y $ are called \textit{partite sets} of $G$, and $ (X, Y) $ is called a bipartition of $ G $. A bipartite graph $ G $ with a bipartition $ (X, Y) $ is a \textit{complete bipartite} if each vertex in $ X $ is adjacent to every vertex in $ Y $ and vice versa.  A complete bipartite graph with a bipartition $ (X, Y) $, where $ |X| = m $ and $ |Y| = n $ is denoted by $K_{m, n}$ (here $|X|$ represents the number of elements in $X$).  Two graphs $G_1$ and $G_2$ are isomorphic, written as $G_1\cong G_2$, if there exists a bijection $f: V(G_1)\longrightarrow V(G_2)$ such that $(u, v) \in E(G_1)$ if and only if $(f(u), f(v)) \in E(G_2)$.

The \textit{union} of $ G_1 $ and $ G_2 $, denoted by $G_1\cup G_2$ is the graph with $V(G_1 \cup G_2) = V(G_1) \cup V(G_2) $ and $E(G_1 \cup G_2) = E(G_1) \cup E(G_2) $. The \textit{join} of $ G_1 $ and $ G_2 $, denoted by $G_1\vee G_2$ is the graph with $ V(G_1 \vee G_2) = V(G_1) \cup V(G_2)$ and $E(G_1 \vee G_2) = E(G_1 \cup G_2)\, \cup\, \{(x, y): x \in V(G_1)\, \text{and}\, y \in V(G_2)\} $.

We now discuss some basic definitions and results for frames in $\mathbb R^n$ with the Euclidean inner product, which are useful for the further discussion. A finite sequence of vectors, $\phi=\{f_i\}_{i=1}^m$ $\subseteq \mathbb{R}^n$ $(m\geq n)$ is called a \textit{frame} in $\mathbb{R}^n$ if there exist positive constants $ A $ and $ B $ with $A\leq B<\infty$ such that
\begin{equation}
	A\|f\|^2\leq \sum\limits_{i=1}^m|\langle f,f_i\rangle|^2\leq B\|f\|^2,\,\,\,\, \text{for all}\,\, f\in \mathbb{R}^n.
\end{equation}
The constants $A$ and $B$ are called \textit{lower} and \textit{upper frame bounds}. If $A = B$, then $\phi$ is called a \textit{tight frame} or an \textit{$A$-tight frame}. If $A = B =1$, then $\phi$ is called a \textit{Parseval frame}. If $\phi=\{f_i\}_{i=1}^m$ is an $A$-tight frame in $\mathbb{R}^n$, then $\phi^\prime=\{\frac{f_i}{\sqrt{A}}\}_{i=1}^m$ forms a Parseval frame in $\mathbb{R}^n$. Given a frame  $\phi=\{f_i\}_{i=1}^m$ in $\mathbb{R}^n$, the corresponding \textit{analysis operator}, denoted by $F_\phi$ is the $m\times n$  matrix whose $i^{th}$ row is $f_i^t$ (the transpose of $f_i$). The operator $F_\phi^t$ is called the \textit{synthesis operator}. The \textit{Grammian operator}, $G_\phi$ is defined as $G_\phi=F_\phi F^t_\phi$ and the \textit{frame operator}, $S_\phi$ is defined as $ S_\phi =F^t_\phi F_\phi$. The frame operator $S_\phi$ is a positive symmetric invertible operator and it satisfies \begin{equation}
S_\phi(f) =\sum\limits_{i=1}^{m}\langle f, f_i\rangle f_i,\,\, \text{for all}\,\, f\in \mathbb{R}^n. 
\end{equation}
A frame $\phi$ is a \textit{Parseval} (\textit{tight}) frame if and only if $S_\phi=I$ ($S_\phi=\lambda I$, for some $\lambda> 0$). Here $I$ is the identity matrix. A frame $\phi=\{f_i\}_{i=1}^m$ is said to be a \textit{scalable frame} in $\mathbb{R}^n$ if there exist non-negative constants $\{a_i\}_{i=1}^m$ for which $\{a_if_i\}_{i=1}^m$ is a Parseval frame. If all the $a_i\text{'s}$ are positive, then the frame $\phi$ is called a \textit{strictly scalable frame} in $\mathbb{R}^n$.

\begin{definition}\cite{Abd18}\label{definition}
    Let $\phi$ be a frame in $\mathbb R^n$. We associate a simple graph $G(\phi)$ whose vertices are the elements of $\phi$ and two distinct vertices $u$ and $v$ are adjacent if and only if  $\langle a, b\rangle\neq 0$.      A simple graph $G$ with $V(G)=m\geq n$ is called a \textit{frame graph} in $\mathbb{R}^n$ if there exists a frame $\phi$ consisting of $m$ vectors in $\mathbb{R}^n$ such that $G(\phi) \cong G$. A frame graph is called a \textit{tight frame graph} in $ \mathbb{R}^n$ if the associated frame is a tight frame in $\mathbb{R}^n$.
\end{definition} 
\begin{definition}\cite{zche14}
     The \textit{representation number} of a given graph $G$ with $m$ vertices is the smallest positive number $n$ such that there exists a frame $\phi=\{f_i\}_{i=1}^m$ in $\mathbb{R}^n$ and $G(\phi)\cong G$. 
\end{definition}

Throughout this paper, we consider only a finite sequence of frame vectors $\phi=\{f_i\}_{i=1}^m$ in $\mathbb R^n$ and $deg(f_i)\geq 1$, for all $f_i\in V(G(\phi))$, which means that  $\phi$ does not contain the zero vector and every member in $\phi$ is not orthogonal to at least one member in $\phi$.  We end the section with some linear algebra results.

 \begin{theorem}\cite{han07} \label{2.1}
Let  $\{f_i\}_{i=1}^m$ be a collection of vectors in $\mathbb R^n$.  Then $\{f_i\}_{i=1}^m$  is a Parseval frame in $\mathbb{R}^n$ if and only if the associated Grammian operator $G_\phi$ is an orthogonal projection of rank $n$.
 \end{theorem}
 \begin{theorem}\cite{hor12}\label{2.2}
 	Let $A$ be an $n\times n$ symmetric matrix. Then the following conditions are equivalent :
 	\begin{enumerate}
 		\item $A$ is positive semi-definite;
 		\item All the eigenvalues of $A$ are non-negative;
 		\item A has decomposition of the form $A=M^tM$, where $M$ is an $m\times n$ matrix with the same rank as $A$.
 	\end{enumerate}
 \end{theorem}
 \begin{theorem}\cite{hor12}\label{2.3}
 	Let $A$ be an $n\times n$ symmetric matrix.  If $\lambda$ is the only eigenvalue of $A$, then $A=\lambda I$.
 \end{theorem}
 \begin{theorem}\cite{hor12}\label{2.4}
 	Let $A$ be an $m\times n$ matrix and $B$ be an $n\times m$ matrix. If $\lambda$ is a non-zero eigenvalue of $AB$ with multiplicity $k$, then $\lambda$ is the eigenvalue of $BA$ with the same multiplicity $k$. 
 \end{theorem}
 
\section{Relationships between graphs and frames}\label{sec3}
In this section, we present some relationships between graphs and frames. Also, we provide some necessary conditions for a frame $\phi$ in $\mathbb{R}^n$ to be a Parseval frame by using the associated graph $G(\phi)$ and the Grammian operator $G_\phi$. 
\begin{proposition}\label{connected}
    Let $\phi=\{f_i\}_{i=1}^m$ be a  frame in $\mathbb{R}^n$ such that $G(\phi)$ is connected. If $\phi^\prime=\{g_i\}_{i=1}^k$  be any collection of non-zero vectors in $\mathbb{R}^n$, then $G(\phi \cup \phi^\prime)$ is a connected graph.
    \begin{proof}
   Assume that $\phi=\{f_i\}_{i=1}^m$ is a frame in $\mathbb{R}^n$ such that $G(\phi)$ is connected. Since $span(\phi)=\mathbb{R}^n$, any non-zero $g\in \mathbb{R}^n$ has $\langle g, f_i\rangle\neq 0$, for some $i\in \{1,2,\ldots,m\}$. Therefore every  vector in $\phi^\prime$ is adjacent to some frame vector in $\phi$. Thus $G(\phi\cup \phi^\prime)$ forms a connected graph.
    \end{proof}
\end{proposition} 

The converse of the above theorem is not true, which is illustrated in the following example.
\begin{example}
    Let $\phi=\{e_i\}_{i=1}^n$ be the standard orthonormal basis in $\mathbb{R}^n$ and let $\phi^\prime=\{e_1+e_2+\cdots+e_n\}$. Clearly, $G(\phi\cup \phi^\prime)$ is a connected graph. But $G(\phi)$ is the totally disconnected graph $\overline{K_n}$. 
\end{example}
\begin{remark}\label{disconnected}
\begin{enumerate}
   \item If $\phi=\{f_i\}_{i=1}^n$ is a tight frame in $\mathbb{R}^n$, then $G(\phi)$ is isomorphic to the totally disconnected graph $\overline{K_n}$. 
   \item If $\phi=\{f_i\}_{i=1}^m$ $(m\geq n)$ is an $A$-tight frame in $\mathbb{R}^n$, then $\phi^\prime=\{\frac{f_i}{\sqrt{A}}\}_{i=1}^m$ forms a Parseval frame in $\mathbb{R}^n$ and $G(\phi^\prime)$ is isomorphic to $G(\phi)$.
   \end{enumerate}
\end{remark}
\begin{theorem}\label{kn}
   Let $\phi=\{f_i\}_{i=1}^n$ be a Parseval frame in $\mathbb{R}^{n-1}$ with the associated graph $G(\phi)$. Then $\langle f_i, f_j\rangle\neq 0$, for all $i\neq j$ and $G(\phi)$ is isomorphic to $K_n$.
    \begin{proof}
      Assume that $\phi=\{f_i\}_{i=1}^n$ is a Parseval frame in $\mathbb{R}^{n-1}$. Then by Theorem \ref{2.1}, $G_{\phi}$ has the eigenvalue $1$ with multiplicity $n-1$ and the eigenvalue $0$ with multiplicity $1$.  Hence  $I-G_{\phi}$ has the eigenvalue $1$ with multiplicity $1$ and the eigenvalue $0$ with multiplicity $n-1$. Clearly, the matrix $I-G_\phi$ preserves the same zero and non-zero off-diagonal entries corresponding to $G_\phi$. By Theorem \ref{2.2}, $I-G_{\phi}$ can be written as $I-G_{\phi}=xx^t$, where $x$ is a column vector in $\mathbb{R}^n$. If one of the coordinates of $x$ is equal to zero, then $G(\phi)$ contains a vertex $f_i$ with $deg(f_i)=0$. But this is a contradiction to $deg(f_i)\geq 1$, for all $f_i\in V(G(\phi))$. Hence, from the matrix $I-G_\phi$, $\langle f_i, f_j\rangle\neq 0$, for all $i\neq j$.  Thus $G(\phi)$ is isomorphic to $K_n$.  
    \end{proof}
\end{theorem}
\begin{corollary}\label{tkn}
      Let $\phi=\{f_i\}_{i=1}^n$ be a tight frame in $\mathbb{R}^{n-1}$ with the associated graph $G(\phi)$. Then $\langle f_i, f_j\rangle\neq 0$, for all $i\neq j$ and $G(\phi)$ is isomorphic to $K_n$.
     \begin{proof}
         The proof follows from Remark \ref{disconnected} and Theorem \ref{kn}.
     \end{proof}
\end{corollary}
\begin{theorem}\label{compliment}
    Let $\phi=\{f_i\}_{i=1}^m$ $(m>n)$ be a Parseval frame in $\mathbb{R}^{n}$ with the associated graph $G(\phi)$. Then there exists a Parseval frame $\phi^{\prime}=\{g_i\}_{i=1}^m$ in $\mathbb{R}^{m-n}$ such that $G(\phi)$ is isomorphic to $G(\phi^{\prime})$.
    \begin{proof}
        Assume that $\phi=\{f_i\}_{i=1}^m$ is a Parseval frame in $\mathbb{R}^n$. Then by Theorem \ref{2.1},  $I-G_{\phi}$ has the eigenvalue $1$ with multiplicity $m-n$ and the eigenvalue $0$ with multiplicity $n$. Clearly, the matrix $I-G_\phi$ preserves the same zero and non-zero off-diagonal entries corresponding to $G_\phi$. Now by Theorem \ref{2.2}, $I-G_{\phi}$ can be written as $I-G_{\phi}=M^tM$, where $M$ is an $(m-n)\times m$ matrix.  The columns of $M$ are denoted by $\phi^{\prime}=\{g_i\}_{i=1}^m$. Hence by Theorems \ref{2.3} and \ref{2.4}, $S_{\phi^\prime}=MM^t=I$. Therefore  $\phi^\prime$ forms a Parseval frame in $\mathbb{R}^{m-n}$ and $G(\phi)$ is isomorphic to $ G(\phi^\prime)$.
    \end{proof}
\end{theorem}
\begin{corollary}\label{tcompliment}
     Let $\phi=\{f_i\}_{i=1}^m$ $(m>n)$ be a tight frame in $\mathbb{R}^{n}$ with the associated graph $G(\phi)$. Then there exists a Parseval frame $\phi^{\prime}=\{g_i\}_{i=1}^m$ in $\mathbb{R}^{m-n}$ such that $G(\phi)$ is isomorphic to $G(\phi^{\prime})$.
     \begin{proof}
         Assume that $\phi$ is a tight frame in $\mathbb{R}^n$ with the frame bound $A$. Then by Remark \ref{disconnected}, $\phi_1=\{\frac{f_i}{\sqrt{A}}\}_{i=1}^m$ forms a Parseval frame in $\mathbb{R}^n$ and $G(\phi)$ is isomorphic to $G(\phi_1)$. Therefore, by Theorem \ref{compliment}, there exists a Parseval frame $\phi^\prime=\{g_i\}_{i=1}^m$ such that $G(\phi_1)$ is isomorphic to $G(\phi^\prime)$ and hence, from the transitive property of isomorphism, $G(\phi)$ is isomorphic to $ G(\phi^\prime)$.
     \end{proof}
\end{corollary}
\begin{theorem}\label{alp}
     Let $\phi=\{f_i\}_{i=1}^n$ be a Parseval frame in $\mathbb{R}^{n-k}$ $(0<k<n)$ with the associated graph $G(\phi)$. Then $\alpha(G(\phi))\leq k$.
    \begin{proof}
   Assume that  $\phi=\{f_i\}_{i=1}^n$ is a tight frame in $\mathbb{R}^{n-k}$. Then by Theorem \ref{compliment}, there exists a tight frame $\phi^\prime=\{g_i\}_{i=1}^n$ in $\mathbb{R}^k$ such that $G(\phi)$ is isomorphic to $G(\phi^\prime)$ and $\alpha(G(\phi))=\alpha(G(\phi^\prime))$. Suppose that  $\alpha(G(\phi))> k$.  Then $\phi^\prime$ contains at least $k+1$ non-zero linearly independent orthogonal vectors in $\mathbb{R}^k$. But this is not possible because the dimension of $\mathbb{R}^k$ is $k$. Therefore $\alpha(G(\phi))\leq k$.
        \end{proof}
\end{theorem}
\begin{corollary}\label{orthogonal}
    Let $\phi=\{f_i\}_{i=1}^n$ $(n\geq 3)$ be a frame in $\mathbb{R}^{n-k}$ $(0<k<n)$. If $\phi$ contains $k+1$ orthogonal vectors in $\mathbb{R}^{n-k}$, then $\phi$ is not a Parseval frame in $\mathbb{R}^{n-k}$.
    \begin{proof}
    Assume that the $k+1$ orthogonal vectors form an orthogonal set in $\mathbb{R}^{n-k}$. Then the independent number of $\alpha(G(\phi))$ is at least $k$. Therefore by Theorem \ref{alp}, $\phi$ is not a Parseval frame in $\mathbb{R}^{n-k}$.
\end{proof}
\end{corollary}
\begin{theorem}\label{diam}
  Let $\phi=\{f_i\}_{i=1}^n$ $(n\geq 3)$ be a Parseval frame in $\mathbb{R}^{n-2}$ with the associated graph $G(\phi)$ is connected. Then $diam(G(\phi))\leq 2$.  
  \begin{proof}
    Assume that $\phi=\{f_i\}_{i=1}^n$ is a Parseval frame in $\mathbb{R}^{n-2}$. Then by Theorem \ref{compliment}, there exists a Parseval frame $\phi^\prime=\{g_i\}_{i=1}^n$ in $\mathbb{R}^2$ such that $G(\phi)$ is isomorphic $G(\phi^\prime)$ and $diam(G(\phi))=diam(G(\phi^\prime))$. Suppose that  $diam(G(\phi^\prime))\geq 3$.  Then without loss of generality, there exists an induced path, $P=\{g_1, g_2,\ldots, g_k\}$ in $G(\phi^\prime)$, with $k\geq 4$. Now, choose the first four vertices from $P$. 
    
        Clearly, the vertices $g_1$ and $g_3$ are not adjacent in $G(\phi^\prime)$, and hence  $g_1$ and $g_3$ form an orthogonal basis in $\mathbb{R}^2$. Therefore, the vertex $g_2$  must be a linear combination of both $g_1$ and $g_3$.   But this is not possible because the vertices $g_2$ and $g_4$ are not adjacent in $G(\phi^\prime)$. Therefore $diam(G(\phi)^\prime)\leq 2$, and hence $diam(G(\phi))\leq 2$. 
  \end{proof}
\end{theorem}
\begin{theorem}\label{bipartite}
 Let $\phi=\{f_i\}_{i=1}^m$ $(m\geq n)$ be a Parseval frame in $\mathbb{R}^{n}$ with associated graph $G(\phi)$. If $G(\phi)$ is a bipartite graph with a bipartition $(X, Y)$, then $|X|=|Y|$.
    \begin{proof}
        It is given that $\phi=\{f_i\}_{i=1}^m$ is a Parseval frame in $\mathbb{R}^n$ and $G(\phi)$ is a bipartite graph with a bipartition $(X, Y)$, where $|X|=k_1$,  $|Y|=k_2$ and $k_1+k_2=m$.  Without loss of generality, we assume that the rows and columns of $G_\phi$ are indexed by both $X$ and $Y$.  
        
        Let us consider the matrix $G_{\phi}=\begin{bmatrix}
            D_1 & B\\
            B^t & D_2
        \end{bmatrix}$. Here $D_1$ is a $k_1\times k_1$ diagonal matrix with rows and columns are indexed by $X$, $D_2$ is a $k_2\times k_2$ diagonal matrix with rows and columns are indexed by $Y$ and $B$ is a $k_1\times k_2$ matrix with rows and columns are indexed by $X$ and $Y$ respectively. Since $\phi$ is a Parseval frame, by Theorem \ref{2.1}, $B^tB$ and $BB^t$ are diagonal. Therefore, the rows and the columns are pairwise orthogonal in $B$ and hence $|X|=|Y|$.
    \end{proof}
\end{theorem}
\begin{corollary}\label{tbipartite}
 Let $\phi=\{f_i\}_{i=1}^m$ $(m\geq  n)$ be a tight frame in $\mathbb{R}^n$ with the associated graph $G(\phi)$. Then the following statements hold good :
 \begin{enumerate}
     \item If $G(\phi)$ is a bipartite graph with a  bipartition $(X, Y)$, with $|X|=k_1$ and $|Y|=k_2$, then $|X|=|Y|$.
     \item If $\phi$ contains an odd number of frame vectors, then the associated graph $G(\phi)$ cannot be a bipartite graph.
 \end{enumerate}
  \begin{proof}
      The proof  follows from Remark \ref{disconnected} and Theorem $\ref{bipartite}$.
  \end{proof}
\end{corollary}
\begin{theorem}\label{independent}
Let $\phi=\{f_i\}_{i=1}^m$ $(m\geq n)$ be a Parseval frame in $\mathbb{R}^n$ with the associated graph $G(\phi)$. If $X$ is an independent set of vertices in $G(\phi)$, then $|X|\leq \frac{m}{2}$.
     \begin{proof}
        Assume that $\phi=\{f_i\}_{i=1}^m$ is a Parseval frame and $X$ is an independent set in $G(\phi)$. Let
         $G_{\phi}=\begin{bmatrix}
            D & B\\
            B^t & A
        \end{bmatrix}$, where $D$ is a diagonal matrix with rows and columns are indexed by $X$, $B$ is a matrix whose rows and columns are indexed by $X$ and $X^c$ respectively, and $A$ is a square matrix with rows and columns are indexed by $X^c$. Since $\phi$ is a Parseval frame, by Theorem \ref{2.1}, $BB^t=D-D^2$. Therefore, the rows are pairwise orthogonal in $B$, so  $|X|\leq |X^c|$. Hence $|X|\leq \frac{m}{2}$.
     \end{proof}
\end{theorem}
\section{Necessary conditions for scalable frames}\label{sec4}
In this section, we give some necessary conditions for scalability of a given frame $\phi=\{f_i\}_{i=1}^m$ in $\mathbb{R}^n$ using its associated graph $G(\phi)$.
\begin{theorem}\label{iso1}
    Let $\phi=\{f_i\}_{i=1}^m$ $(m\geq n)$ be a frame in $\mathbb{R}^n$ with the associated graph $G(\phi)$. If $\phi$ is a scalable frame in $\mathbb{R}^n$, then the graph corresponding to that scalable frame is a subgraph of $G(\phi)$.
    \begin{proof}
        Assume that $\phi=\{f_i\}_{i=1}^m$ is a scalable frame in $\mathbb{R}^n$. Then there exist non-negative constants $\{a_i\}_{i=1}^m$, such that $\{a_if_i\}_{i=1}^m$ forms a tight frame. After dropping those vectors that are scaled with the zero constant, the vectors $\{a_if_i\}_{i=k+1}^m$ preserve the adjacency relation between the vectors $\{f_i\}_{i=k+1}^m$. Therefore, the scalable frame corresponding to $\phi$ forms a subgraph of $G(\phi)$.
    \end{proof}
\end{theorem}
\begin{corollary}\label{iso2}
     Let $\phi=\{f_i\}_{i=1}^m$ $(m\geq n)$ be a frame in $\mathbb{R}^n$ with the associated graph $G(\phi)$. If $\phi$ is a strictly scalable frame in $\mathbb{R}^n$, then the graph corresponding to that strictly scalable frame is isomorphic to $G(\phi)$.
    \end{corollary}

\begin{theorem}\label{snecessary}
  Let $\phi=\{f_i\}_{i=1}^m$ $(m\geq 3)$ be a frame in $\mathbb{R}^n$ and the associated graph $G(\phi)$ be a connected graph. If $G(\phi)$ contains two non-adjacent vertices that share exactly one common neighbor, then $\phi$ is not strictly scalable frame in $\mathbb{R}^n$.
  \begin{proof}
      Suppose that $\phi=\{f_i\}_{i=1}^m$ is a strictly scalable frame in $\mathbb{R}^n$. Then there exist positive constants $\{a_i\}_{i=1}^m$ such that $\phi^\prime=\{a_if_i\}_{i=1}^m$ forms a Parseval frame in $\mathbb{R}^n$. Let $f_m$ and $f_n$ be  two non-adjacent vertices in $G(\phi)$ with the common neighbor $f_\ell$. By Corollary \ref{iso2}, $a_mf_m$ and $a_nf_n$ are also two non-adjacent  vertices in $G(\phi^\prime)$ with the common neighbor $a_\ell f_\ell$. The vector $a_mf_m$ can be written in the following way :
     \begin{equation*}
         a_mf_m=\sum\limits_{i=1}^m\langle a_mf_m, a_if_i\rangle a_if_i \,\,\,\text{and}\,\,\, \langle a_mf_m, a_nf_n\rangle=\langle a_mf_m, a_\ell f_\ell\rangle \langle a_\ell f_\ell, a_nf_n\rangle.
     \end{equation*}
    Since $\langle a_mf_m, a_nf_n\rangle=0$, either the value of $\langle a_mf_m, a_\ell f_\ell\rangle$ is zero or the value of $\langle a_\ell f_\ell, a_nf_n\rangle$ is zero. But this is not possible, because the vertices $\{a_mf_m, a_\ell f_\ell\}$ and $\{a_\ell f_\ell, a_nf_n\}$ are adjacent in $G(\phi^\prime)$. Hence $\phi=\{f_i\}_{i=1}^m$ is not a strictly scalable frame in $\mathbb{R}^n$.
  \end{proof}
  \end{theorem}
  \begin{corollary}
       Let $\phi=\{f_i\}_{i=1}^m$ $(m\geq 3)$ be a frame in $\mathbb{R}^n$ 
       and the associated graph $G(\phi)$ be a connected graph. If $G(\phi)$ contains either a leaf vertex or a bridge, then $\phi$ is not a strictly scalable frame in $\mathbb{R}^n$.
       \begin{proof}
           Assume that $\phi=\{f_i\}_{i=1}^m$ is a strictly scalable frame in $\mathbb{R}^n$. Then there exist positive constants $\{a_i\}_{i=1}^m$ such that $\phi^\prime=\{a_if_i\}_{i=1}^m$ forms a Parseval frame in $\mathbb{R}^n$. Let $f_m$ be a leaf vertex in $G(\phi)$ with the neighbor $f_n$. By Corollary \ref{iso2}, $a_mf_m$ is a leaf vertex in $G(\phi^\prime)$ with the neighbor $a_nf_n$. Since $G(\phi^\prime)$ is a connected graph with at least $3$ vertices,  there exists a vertex $a_\ell f_\ell$ in $G(\phi^\prime)$ such that $a_\ell f_\ell$ is adjacent to $a_nf_n$. Hence by Theorem \ref{snecessary}, $\phi$ is not a strictly scalable frame in $\mathbb{R}^n$. The proof is similar if $G(\phi)$ contains a bridge.
       \end{proof}
   
     \end{corollary}

  \begin{corollary}
      Let $\phi=\{f_i\}_{i=1}^m$ $(m\geq 3)$ be a frame in $\mathbb{R}^n$ with the disconnected associated graph $G(\phi)$ . If one of the components of $G(\phi)$ contains at least $3$ vertices and contains either leaf or bridge, then $\phi$ is not a strictly scalable frame in $\mathbb{R}^n$.
  \end{corollary}
  \begin{corollary}
       Let $\phi=\{f_i\}_{i=1}^m$ $(m\geq 3)$ be a frame in $\mathbb{R}^n$ with the connected associated graph $G(\phi)$. If $\phi$ is strictly scalable, then all the edges in $G(\phi)$ must lie in the cycle. 
  \end{corollary}
  \begin{theorem}
       Let $\phi=\{f_i\}_{i=1}^m$ $(m\geq 3)$ be a frame in $\mathbb{R}^n$ with the associated graph $G(\phi)$. Then $\phi$ is not a strictly scalable frame in $\mathbb{R}^n$ if one of the following conditions holds good :
       \begin{enumerate}
           \item $G(\phi)$ is isomorphic to a tree.
           \item $G(\phi)$ is isomorphic to a bipartite graph with each partite set having a different size.
           \item $G(\phi)$ contains an independent set $X$ with $|X|>\frac{m}{2}$.
       \end{enumerate}
       \begin{proof}
           The proof follows by Theorems \ref{bipartite}, \ref{independent} and \ref{snecessary}.
       \end{proof}
  \end{theorem}

  \begin{theorem}
      Let $\phi=\{f_i\}_{i=1}^n$ $(n\geq 3)$ be a frame in $\mathbb{R}^{n-2}$ with the associated graph $G(\phi)$. If $\phi$ is a strictly scalable in $\mathbb{R}^{n-2}$, then the following conditions holds good : 
          \begin{enumerate}
              \item $\alpha(G(\phi))\leq 2$.
              \item $diam(G(\phi))\leq 2$, if $G(\phi)$ is a connected graph.
              \item $\phi$ does not contains an orthogonal set with three elements in $\mathbb{R}^{n-2}$.
          \end{enumerate}
          \begin{proof}
              The proof follows by Theorems \ref{alp} and \ref{diam}.
          \end{proof}
  \end{theorem}
  \begin{example}
      Let $\phi=\{f_i\}_{i=1}^5$ be a frame in $\mathbb{R}^\ell$ and $G(\phi)$ be the graph corresponding to $\phi$ as given in Figure \ref{ss10}.
     \begin{figure}[h]
	\centering
 \includegraphics[width=0.35\linewidth]{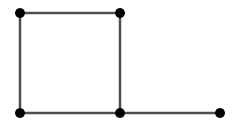}
	\caption[ft1]{$G(\phi)$} 
	\label{ss10}
\end{figure}

Since $\alpha(G(\phi))=3$ and dimension of $span(\phi)$ is at least 3, the representation number of $G(\phi)$ is at least $3$. Hence by Theorem \ref{snecessary}, $\phi$ is not a strictly scalable frame in $\mathbb{R}^\ell$, where $3\leq \ell \leq 5$.
  \end{example}
\begin{remark}
 In \cite{cas23}, the authors present a method to find the scaling constants of a given frame in $\mathbb{R}^n$, by solving a system of linear equations.   Solving   systems of linear equations in higher dimensions is not  easy. 
In this case, given a frame $\phi=\{f_i\}_{i=1}^n$ in $\mathbb{R}^\ell$, we can use our method to verify the scalability. Especially, the value of $\ell$ is $(n-2)$ or $(n-1)$ or $n$. 
In our approach, we can conclude existence of  scaling constants, by using the graph corresponding to $\phi$.

\end{remark}

The following examples are given in \cite{cas23}, and the authors have  verified scalability of frames by solving systems of linear equations. Here, we verify the same by using graphs.
\begin{example}
    Let $M_1$ be a matrix whose columns $\phi_1=\{f_i\}_{i=1}^4$ form a frame in  $\mathbb{R}^4$ and let $G(\phi_1)$ be the graph corresponding to $\phi_1$, as given in Figure \ref{img3}.
    
    \begin{figure}[h]
\centering
\begin{minipage}[b]{.45\textwidth}
\includegraphics[width=\textwidth]{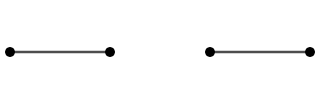}
\caption{$G(\phi_1)$}\label{img3}
\end{minipage}\hfill
\begin{minipage}[b]{.45\textwidth}
\centering
\small
\(
\vspace{.5 cm}
M_1= 
\begin{pmatrix}
        1 & 1 & 0 & 0\\
	  2 & -2 & 0 & 0 \\
		0 & 0 & 1 & 1 \\
        0 & 0 & 2 & -2\\
\end{pmatrix}
\)
\end{minipage}
\end{figure}
   
 Suppose that  $\phi_1$ is a scalable frame. Then there exist positive constants $\{a_i\}_{i=1}^4$ such that $\phi_1^\prime=\{a_if_i\}_{i=1}^4$ forms a tight frame in $\mathbb{R}^4$.  By Theorem \ref{iso2}, $G(\phi_1^\prime)$ is isomorphic to $G(\phi_1)$. But this is not possible by Remark \ref{disconnected}, because $G(\phi_1^\prime)$ is isomorphic to $\overline{K_4}$.
\end{example}
\begin{example}
	
	    Let $M_2$ be a matrix whose columns $\phi_2=\{f_i\}_{i=1}^4$ form a frame in  $\mathbb{R}^4$ and let $G(\phi_2)$ be the graph corresponding to $\phi_2$, as given in Figure \ref{img4}.
   \begin{figure}[ht]
\centering
\begin{minipage}[b]{.45\textwidth}
\hspace{1cm}
\vspace{-.4cm}
\includegraphics[width=0.7\textwidth]{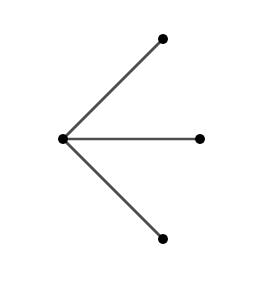}
\caption{$G(\phi_2)$}\label{img4}
\end{minipage}\hfill
\begin{minipage}[b]{.45\textwidth}
\centering
\small
\(
\vspace{.5cm}
M_2= \begin{pmatrix}
    1 & -1 & 1 & 1\\[.5em]
		1 & 1 & -1 & 1 \\[.5em]
		1 & 1 & 1  & -1 \\[.5em]
        1 & 1 & 1 & 1   
\end{pmatrix}
\)
\end{minipage}
\end{figure}

Suppose that  $\phi_2$ is a scalable frame. Then there exist positive constants $\{a_i\}_{i=1}^4$ such that $\phi_2^\prime=\{a_if_i\}_{i=1}^4$ forms a tight frame in $\mathbb{R}^4$.  By Theorem \ref{iso2}, $G(\phi_2^\prime)$ is isomorphic to $G(\phi_2)$. But this is not possible by Theorem \ref{snecessary}.
\end{example}
\begin{example}\label{4.8}
Let $\phi=\phi_1\cup \phi_2$ be a frame for $\mathbb{R}^4$ which is given by the columns of the following matrix, and Figure \ref{img5} be the graph corresponding to $\phi$.
  \begin{figure}[h]
\centering
\begin{minipage}[b]{.45\textwidth}
\hspace{.5cm}
\includegraphics[width=0.7\textwidth]{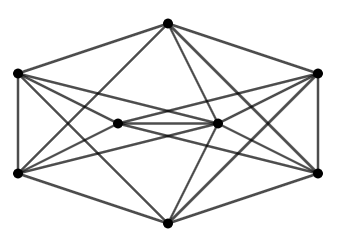}
\caption{$G(\phi)$}\label{img5}
\end{minipage}\hfill
\begin{minipage}[b]{.5\textwidth}
\centering
\small
\label{M}
\(
\vspace{.5cm}
M=\begin{pmatrix}
     1 & 1 & 0 & 0 & 1 & -1 & 1 & 1\\
  2 & -2 & 0 & 0 & 1 & 1 & -1 & 1 \\
  0 & 0 & 1 & 1 & 1 & 1 & 1  & -1 \\
  0 & 0 & 2 & -2 & 1 & 1 & 1 & 1
\end{pmatrix}
\)
\end{minipage}
\end{figure}

The graph corresponding to $\phi$ is isomorphic to $(K_2\cup K_2)\vee K_{1,3}$ and the dimension of $span(\phi)$ is $4$. Clearly, this graph satisfies all the necessary conditions in Theorem $\ref{snecessary}$. 

Assume that  $\phi=\{f_i\}_{i=1}^8$ is a frame in $\mathbb{R}^\ell$, where $\ell\geq 7$ and $G(\phi)$ is isomorphic to Figure $\ref{img5}$.  Then we can conclude that $\phi$ is not a scalable frame in $\mathbb{R}^\ell$ by using graph which is explained in the following section.
\end{example}
\begin{theorem}\label{edgein}
 Let $\phi=\{f_i\}_{i=1}^m$ be a frame for $\mathbb{R}^{n}$ and $G(\phi)$ be the graph corresponding to $\phi$. If $N[f_i]\neq N[f_j]$, for all $f_i\neq f_j\in V(G(\phi))$, then any two adjacent vectors of $\phi$ form an independent set in $\mathbb{R}^n$.
 \begin{proof}
  Assume that  there exist two adjacent vertices $f_i$ and $f_j$ in $V(G(\phi))$, form a linearly dependent set in $\mathbb{R}^n$. This implies that $f_i$  is a linear combination of $f_j$, hence $N[f_i]=N[f_j]$ in $G(\phi)$.
 \end{proof}
 \end{theorem}

 \subsection{Scalability of special classes of frames}\label{sec5}
In this section, we verify the scalability of a given frame by using some special graph classes. The following lemma is mentioned in \cite{zche14}. By using the lemma and the theorems in the previous sections, we can verify the scalability of some particular classes of frames using graphs.
\begin{lemma} \cite{zche14} \label{tree}
Let $\phi=\{f_i\}_{i=1}^n$ be a collection of vectors in $\mathbb{R}^n$. If $G(\phi)$ is isomorphic to a tree with $n$ vertices, then any sub-collection of $n-1$ vectors form an independent set in $\mathbb{R}^n$.
\end{lemma}
\begin{theorem}
    Let $\phi=\{f_i\}_{i=1}^m$ be a frame in $\mathbb{R}^n$ with the associated graph $G(\phi)$. If $\phi$ is a strictly scalable frame in $\mathbb{R}^n$, then the length of the induced path in $G(\phi)$ is not more than $\lfloor \frac{n}{2}\rfloor+1$.
    \begin{proof}
        Assume that $\phi=\{f_i\}_{i=1}^m$ is a strictly scalable frame in $\mathbb{R}^n$. Then there exist  positive constants $\{a_i\}_{i=1}^m$ such that $\phi^\prime=\{a_if_i\}_{i=1}^m$ forms a Parseval frame in $\mathbb{R}^n$. Suppose that $G(\phi)$ contains an induced path, $P=\{f_1, f_2,\ldots, f_{\lfloor \frac{n}{2}\rfloor +1}\}$ of length more than $\lfloor \frac{n}{2}\rfloor +1$. Then by Corollary \ref{iso2},  $G(\phi^\prime)$ also contains an induced path, $P^\prime=\{a_1f_1, a_2f_2,\ldots, a_{\lfloor \frac{n}{2}\rfloor +1}f_{\lfloor \frac{n}{2}\rfloor +1}\}$ of length more than $\lfloor \frac{n}{2}\rfloor +1$ and by Lemma \ref{tree}, the representation number of $G(\phi^\prime)$ is at least $\lfloor \frac{n}{2}\rfloor +1$. But this is not possible by Theorem \ref{compliment}. Hence the length of the induced path in $G(\phi)$ is not more than $\lfloor \frac{n}{2}\rfloor +1$. 
    \end{proof}
\end{theorem}
\begin{theorem}\label{cycle}
	Let $\phi=\{f_i\}_{i=1}^n$ $(n\geq 7)$ be a frame in $\mathbb{R}^\ell$. If $G(\phi)$ is isomorphic to a cycle with $n$ vertices, then $\phi$ is not a scalable frame in $\mathbb{R}^\ell$. 
	\begin{proof}
		Suppose that $\phi$ is scalable in $\mathbb{R}^\ell$. Then there exist non-negative constants $\{a_i\}_{i=1}^{n}$ such that $\phi^\prime=\{a_if_i\}_{i=1}^n$ forms a tight frame in $\mathbb{R}^\ell$ and by Theorem \ref{iso1}, $G(\phi^\prime)$ is isomorphic to subgraph of $G(\phi)$. Since $C_n$ contains an induced path of length $n-2$. Without loss of generality, we assume that the vectors $S=\{f_1, f_2,\ldots, f_{n-1}\}$ form an induced path in $C_n$. Then, by Lemma \ref{tree}, any subcollection of $(n-2)$ vectors from $S$ form an independent set in $\mathbb{R}^\ell$. Therefore the representation number of $G(\phi)$ is at least $(n-2)$. 
		
		If $\ell=n$, then by Remark \ref{disconnected}, $\phi$ is not a strictly scalable frame.
		
		If $\ell=(n-1)$, then by Theorem \ref{kn}, $\phi$ is not a strictly scalable frame. Therefore, exactly one of the scaling constants is zero. Suppose that  there exists non-negative constants $\{a_i\}_{i=1}^n$ such that $\phi^\prime=\{a_if_i\}_{i=1}^n$ forms a tight frame in $\mathbb{R}^{n-1}$. Then the representation number of $G(\phi^\prime)$ is at least $5$. Therefore, by Corollary \ref{tcompliment}, this is not possible. Hence $\phi$ is not a scalable frame in $\mathbb{R}^{n-1}$.  
		
		If $\ell=n-2$, then by Theorem \ref{compliment}, $\phi$ is not a strictly scalable frame. Therefore, either one or two scaling constants must be equal to zero. Suppose that  exactly one of the scaling constants is equal to zero. Then, the graph $G(\phi^\prime)$ contains an induced path of length $(n-1)$. But by Theorem \ref{compliment}, this is not possible. Suppose that  two of the scaling constants are zero, then the representation number of $G(\phi^\prime)$ is at least $3$. Therefore, by Theorem \ref{compliment}, this is also not possible. Hence $\phi=\{f_i\}_{i=1}^n$ is not a scalable frame in $\mathbb{R}^{n-2}$.
	\end{proof}
\end{theorem}

In the following two examples, we illustrate how to verify the scalability of a given frame $\phi$ in $\mathbb{R}^n$ using graphs.
\begin{example}\label{5.2}
	Let $\phi=\{f_i\}_{i=1}^7$ be a frame in $\mathbb{R}^\ell$ and Figure \ref{img7} be the graph corresponding to $\phi$. Then by Theorem \ref{edgein}, the vectors $\{f_1,f_2,f_3,f_4\}$ form an independent set in $\mathbb{R}^\ell$. Therefore the representation number of $G(\phi)$ is at least $4$. Suppose that  $\phi$ is a scalable frame in $\mathbb{R}^\ell$. Then there exist  non-negative constants $\{a_i\}_{i=1}^{7}$ such that $\phi^\prime=\{a_if_i\}_{i=1}^{7}$ forms a tight frame in $\mathbb{R}^\ell$ and by Theorem \ref{iso1}, $G(\phi^\prime)$ is isomorphic to a subgraph of $G(\phi)$.
	\begin{figure}[h]
		\centering
		\includegraphics[width=0.3\linewidth]{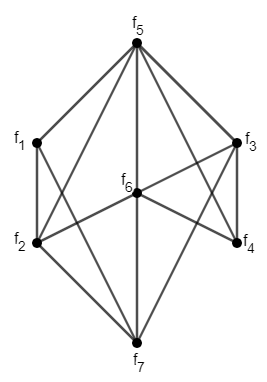}
		\caption[ft1]{$G(\phi)$} 
		\label{img7}
	\end{figure}
	
	If the value of $\ell$ is $7$, then by Remark \ref{disconnected} and Theorem \ref{iso2}, $\phi$ is not a strictly scalable frame in $\mathbb{R}^7$.
	
If the value of $\ell$ is $6$, then by Theorem \ref{kn}, $\phi$ is not a strictly scalable frame. Therefore, if $\phi$ is a scalable frame in $\mathbb{R}^6$, then exactly one of the scaling constants $\{a_i\}_{i=1}^7$ must be  zero. Suppose that  one of the scaling constants corresponding to $\{f_2, f_3, f_5, f_6, f_7\}$ is zero. Then the graph $G(\phi^\prime)$ contains non-adjacent vertices with exactly one common neighbor. But by Theorem \ref{snecessary}, this is not possible. Suppose that  one of the scaling constants corresponding to $\{f_1, f_4\}$ is  zero. Then, by Theorem \ref{edgein}, the representation number of $G(\phi^\prime)$ is at least $3$. But this is not possible by Theorem \ref{compliment}. Hence $\phi$ is not a scalable frame in $\mathbb{R}^6$. 
	
	If the value of $\ell$ is $5$, then by  Theorems  \ref{compliment} and \ref{edgein}, $\phi$ is not a strictly scalable frame in $\mathbb{R}^5$. Therefore, if $\phi$ is scalable in $\mathbb{R}^5$, then either one or two scaling constants $\{a_i\}_{i=1}^7$ must be zero. Suppose that  one of the scaling constants corresponding to $\phi=\{f_i\}_{i=1}^7$ is  zero. Then  the representation number of $G(\phi^\prime)$ is at least $3$. But this is not possible by Theorem \ref{compliment}. Hence if $\phi=\{f_i\}_{i=1}^7$ is a scalable frame in $\mathbb{R}^5$, then precisely two of the scaling constants must be zero.
	
	If the value of $\ell$ is $4$, then by  Theorems \ref{compliment} and  \ref{edgein}, $\phi$ is not a strictly scalable frame in $\mathbb{R}^4$. Therefore, if $\phi=\{f_i\}_{i=1}^7$ is a scalable frame in $\mathbb{R}^4$, then at least one of the scaling constants must be  zero.
\end{example}
\begin{example}
	Let $\phi=\{f_i\}_{i=1}^8$ be a frame in $\mathbb{R}^\ell$ and Figure \ref{img11} be the graph corresponding to $\phi$. Suppose that $\phi$ is a scalable frame in $\mathbb{R}^\ell$. Then there exist  non-negative constants $\{a_i\}_{i=1}^8$ such that $\phi^\prime=\{a_if_i\}_{i=1}^8$ forms a tight frame in $\mathbb{R}^\ell$ and by Theorem \ref{iso1}, $G(\phi^\prime)$ is isomorphic to a subgraph of $G(\phi)$. In the case of $\ell \geq 7$, $\phi$ is not a scalable frame, because the representation number of $G(\phi^\prime)$ at least $3$. Therefore, $\phi=\{f_i\}_{i=1}^8$ is not a scalable frame in $\mathbb{R}^\ell$, where $\ell \geq 7$. Suppose that  $\phi$ is a scalable frame in $\mathbb{R}^6$. Then by Theorem \ref{compliment}, one of the scaling constants corresponding to $\{f_3, f_4, f_5\}$ must be zero.
	\begin{figure}[h]
		\centering \includegraphics[width=0.5\linewidth]{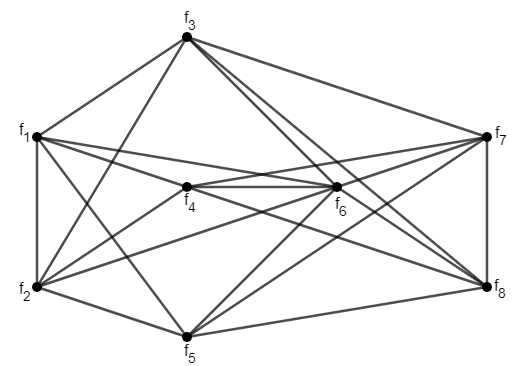}
		\caption[ft3]{$G(\phi)$} 
		\label{img11}
	\end{figure}
\end{example}

\section{Acknowledgements}
We are thankful to Dr. I. Jeyaraman  from National Institute of Technology Tiruchirappali for a number of interesting discussions related to this topic and insights on this paper. The present work of the second author was partially supported by Science and Engineering Research Board (SERB), Department of Science and Technology, Government of India (Reference Number: MTR/2023/000471) under the scheme \enquote{Mathematical Research Impact Centric Support (MATRICS)}.

\end{document}